\documentclass[global]{svjour}  

\usepackage{graphicx} 
\usepackage{amsfonts}
\usepackage{amsmath}
\usepackage{amssymb}

\newtheorem{hypothesis}{Hypothesis}

\newcommand{\at}{{\char '100}}
\newcommand{\EXP}[1]{ \mbox{\large e}^{#1}}

\newcommand{\imat}{{\mathrm{i}}}

\newcommand{\kete}[1]{|\kern.3ex#1\kern.3ex\rangle}

\newcommand{\brae}[1]{\langle\kern.3ex #1 \kern.3ex|} 

\newcommand{\RR}{\mathbb{R}}
\newcommand{\ZZ}{\mathbb{Z}}

\newcommand{\NN}{\mathbb{N}}
\newcommand{\CC}{\mathbb{C}}

\renewcommand{\Re}{\mathrm{Re}}
\renewcommand{\Im}{\mathrm{Im}}

\begin{document}


\title{Upper and lower bounds for an eigenvalue associated with a positive eigenvector}
\author{Amaury Mouchet\\ mouchet\at phys.univ-tours.fr}
\institute{Laboratoire de Math\'ematiques et de Physique Th\'eorique 
                            \textsc{(cnrs umr 6083)}, Universit\'e Fran\c{c}ois Rabelais
                            Ave\-nue Monge, Parc de Grandmont 37200 Tours, 
                            France.}
\date{\today}

\maketitle

\begin{abstract} When an eigenvector of a semi-bounded operator is
  positive, we show that a remarkably simple argument allows to obtain upper and
  lower bounds for its associated eigenvalue. This theorem is a
  substantial generalization of Barta-like inequalities and can be
  applied to non-necessarily purely quadratic Hamiltonians.  An
  application for a magnetic Hamiltonian is given and the
  case of a discrete Schr\"odinger operator is also discussed.  It is
  shown how this approach leads to some explicit bounds on the
  ground-state energy of a system made of an arbitrary number of
  attractive Coulombian particles.
\end{abstract}

\section{Introduction}

In most situations, the principal eigenvalue of a semi-bounded
operator cannot be obtained explicitly whereas it plays a crucial role
in physics: the smallest vibration frequency of an elastic system, the
fundamental mode of an electromagnetic cavity, the ground-state energy
of a quantum system with a finite number of degrees of freedom, the
energy of the vacuum in a quantum field theory, the equilibrium state
at zero temperature in statistical physics, etc.  There are actually
very few ways --- which are usually specific to a restricted class of
systems\footnote{Finding lower bounds for the smallest eigenvalue of a
  typical Hamiltonian is far much difficult than finding upper bounds.
  For successful attempts, see for instance the moment method proposed
  in \cite{Handy/Bessis85a}, the Riccati-Pad\'e method proposed in
  \cite{Fernandez+89a} and the lower bounds obtained for few-body
  systems in \cite{Benslama+98a}.} --- to obtain accurate
approximations of an eigenvalue with a rigorous control on the errors
and a reasonable amount of numerical computations. For instance, in a
typical Dirichlet-Laplacian problem defined for an open connected
set~$\mathcal{Q}\subset\mathbb{R}^d$; $d\in\NN$; Barta's
inequalities~\cite{Barta37a} allow to bound to the lowest
eigenvalue~$e_0$: The determination of a lower (resp. upper) bound
requires the finding of the absolute minimum (resp. maximum) of a smooth
function defined on~$\mathcal{Q}$.  Compared to the general and
traditional methods like the Rayleigh-Schr\"odinger perturbative
series and the Rayleigh-Ritz or Temple variational methods, an
advantage of Barta's approach is not only to naturally provide both an
upper and a lower bound, but also does not involve the calculation of
any integral. Therefore, generalizations of Barta's inequalities can
lead to interesting spectral information. This generalization
 has been carried out
in two directions:

(i) for Laplacian operators acting on square integrable
functions defined on a Riemannian manifold (for a recent work on this 
subject, see \cite{Pacelli/Montenegro04a});

 (ii) for Schr\"odinger 
operators of the form~$-\Delta+V$ acting on square integrable
functions defined on an open set of~$\mathbb{R}^d$
\cite{Barnsley78a,Baumgartner79a,Thirring79b,Crandall/Reno82a,Schmutz85a}
or more generally second order elliptic operators \cite{Protter/Weinberger66a,Berestycki+94a,Harrell05a}.

In both cases, the proofs of Barta's inequalities involve the use of a
Kato-like inequality and therefore rely extensively on the somehow
specific properties of the purely quadratic differential operator.  In
this paper, I propose a significant extension of Barta's inequalities that will
rely on the properties of one eigenvector only.  More precisely,
with a remarkably simple argument, I
will show that we can obtain upper and lower bounds for the eigenvalue
$e_0$ associated with the eigenvector~$\Phi_0$ of an operator, under
the only hypothesis that~$\Phi_0$ is real and non-negative.  This
result includes cases (i) and (ii) because the Krein-Rutman theorem
guarantees the positivity of~$\Phi_0$ \cite[\S XIII.12]{Reed/Simon78a}
for the smallest eigenvalue~$e_0$\footnote{The positivity of~$\Phi_0$
  is also required for the traditional proofs of Barta's inequalities;
  this explains why, in case (i) and (ii), they concern the lowest
  eigenvalue only.}, but also applies for a much wider class of
operators including

 (iii) the Schr\"odinger
operators involving a magnetic field, \textit{e.g.} the Hydrogen atom in a 
Zeeman configuration; 

(iv) discrete Hamiltonians  \textit{e.g.} 
 the one occurring in the Harper model;

(v) some integral or pseudo-differentiable  
operators, \textit{e.g.}  the Klein-Gordon or spinless Salpeter 
Hamiltonians.

The next section fixes the notations, proves the main general results
(theorems~\ref{th:main} and \ref{th:generalnonsym}).
Section~\ref{sec:oldnew} shows how the original argument given in
section~\ref{sec:mainresult} actually embraces and generalizes the
Barta-like inequalities that have been already obtained in the
literature and furnishes guidelines for numerically improve the bounds
on~$e_0$. Sections~\ref{sec:manybody} and~\ref{sec:discrete} provide
two applications in the differentiable case (many-body problem) and in
the discrete case respectively.

\section{Bounding the principal eigenvalue with the local energy}\label{sec:mainresult}

\subsection{General inequalities}\label{subsec:genineq}

In the following, $\mathcal{Q}$ will be a locally compact space
endowed with a positive Radon measure~$\mu$.
$\langle\psi|\varphi\rangle$ will denote the scalar product between
two elements~$\psi$ and $\varphi$ belonging to the Hilbert space of
the square integrable complex
functions~$\mathsf{L}^2(\mathcal{Q},\mu)$:

$$\langle\psi|\varphi\rangle=\int_\mathcal{Q}\overline{\psi}(q)\,\varphi(q)\;d\mu(q)\;.$$

$\mathsf{D}(H)$ will denote the domain of the operator $H$ acting on~$\mathsf{L}^2(\mathcal{Q},d\mu)$. The crucial
hypothesis on~$H$ is the following:

\begin{hypothesis}\label{hyp:sympos}
  The operator~$H$ is symmetric and has one real
  eigenvector~$\Phi_0\in\mathsf{D}(H)$ such that~$\Phi_0\geqslant0$
  (almost everywhere with respect to $\mu$) on~$\mathcal{Q}$.
\end{hypothesis}

If~$e_0$ stands for the eigenvalue of~$H$ associated with~$\Phi_0$,
the symmetry of~$H$ implies that, for all $\varphi\in\mathsf{D}(H)$, we
have $\langle\Phi_0|(H-e_0)\varphi\rangle=0$ that is

$$
 \forall\varphi\in\mathsf{D}(H), \qquad 
 \int_\mathcal{Q}\overline{\Phi_0}(q)\,(H-e_0)\varphi(q)\;d\mu(q)=0\;.
$$     

Taking the real part of the integral, we can see that the support of
$q\mapsto\Re\big[\overline{\Phi_0}(q)\,(H-e_0)\varphi(q)\big]$ either
is empty, either contains two disjoints open sets~$\mathcal{Q}_\pm$
such that $\Re\big(\overline{\Phi_0}(H-e_0)\varphi\big)\gtrless0$ and
$\mu(\mathcal{Q}_\pm)>0$. The hypothesis of positivity of~$\Phi_0$
implies that on $\mathcal{Q}_\pm$, we have
$\Re\big[(H-e_0)\varphi\big]\gtrless0$. The last results motivates the
following definition:

\begin{definition}[Local energy]
\label{def:localenergy}
 For any $\varphi$ in $\mathsf{D}(H)$, the 
 \emph{local energy} is the function~$E_\varphi:
\mathcal{Q}\to\overline{\mathbb{R}}$ defined by 
\begin{equation}\label{eq:localenergy}
   E_\varphi(q)=\frac{\Re{\big(H\varphi(q)\big)}}{\Re\big(\varphi(q)\big)}\;.
\end{equation}
\end{definition}

Therefore from what precedes, we have obtained the main theorem:

\begin{theorem}\label{th:main}
  For any symmetric operator~$H$ on~$\mathsf{L}^2(\mathcal{Q},\mu)$
  having an eigenvalue~$e_0$ whose corresponding eigenfunction is
  nonnegative almost everywhere on~$\mathcal{Q}$, we have
\begin{equation}\label{eq:inequalities}
\forall\varphi\in\mathsf{D}(H)\ \mathrm{such\ that}\ \Re(\varphi)\geqslant0,\qquad\inf_\mathcal{Q}(E_\varphi)\leqslant e_0 \leqslant\sup_\mathcal{Q}(E_\varphi)\;.
\end{equation}
\end{theorem}

Actually, for a non-symmetric operator~$K$, we can 
 keep working with its adjoint $K^\ast$ and easily 
generalize the
above argument:
 \begin{theorem}\label{th:generalnonsym} Let $K$ being an operator
 on~$\mathsf{L}^2(\mathcal{Q},\mu)$ having 
an eigenvalue~$k_0$ whose corresponding eigenfunction
 is real and nonnegative almost everywhere on~$\mathcal{Q}$, we have
$\forall\varphi\in\mathsf{D}(K^\ast)\ \mathrm{such\ that}\ \varphi>0$,
\begin{subequations}\label{eq:inequalitiesnonsym}
\begin{equation}\label{eq:inequalitiesnonsymRe}
\inf_\mathcal{Q}\left[\frac{\Re(K^\ast\varphi)}{\varphi}\right]\leqslant\ \Re(k_0)\ \leqslant\sup_\mathcal{Q}\left[\frac{\Re(K^\ast\varphi)}{\varphi}\right]\;;
\end{equation}
\begin{equation}\label{eq:inequalitiesnonsymIm}
\inf_\mathcal{Q}\left[-\frac{\Im(K^\ast\varphi)}{\varphi}\right]\leqslant\  \Im(k_0)\  \leqslant\sup_\mathcal{Q}\left[-\frac{\Im(K^\ast\varphi)}{\varphi}\right]\;.
\end{equation}
\end{subequations}

\end{theorem}
This generalization may be of physical relevance. There are some
models (\textit{e.g.} the so-called ``kicked'' systems, or quantized
maps) where the dynamics is described ``stroboscobically''
\textit{i.e.}  implemented by a unitary operator (the Floquet
evolution operator) that cannot be constructed from a smooth
Hamiltonian. However, we will not consider this possibility here, and
up to the end of this paper, $H$ will denote a symmetric operator.

\subsection{Optimization strategy}\label{subsec:strategy}

Since generally, the eigenfunction $\Phi_0$ is not known exactly, it
will be approximated with the help of test functions that belong to a
trial space $\mathsf{T}(H)\subset\mathsf{D}(H)$, very much like the
variational method. Since we want a test function to mimic $\Phi_0$ at
best, we will restrict $\mathsf{T}(H)$ to functions that respect the
\textit{a priori} known properties of $\Phi_0$: its positivity, its
boundary conditions and its symmetries if there are any.  For each
test function the error on $e_0$ is controlled by
inequalities~\eqref{eq:inequalities}.  Therefore, the strategy for
obtaining reasonable approximations is clear: First, we must choose or
construct $\varphi$ to eliminate all the singularities of the local
energy in order to work with a bounded function and second, perturb
the test function in the neighborhood of the absolute minimum (resp.
maximum) of the local energy in order to increase (resp. decrease) its
value. For practical and numerical computations, this perturbation
will be implemented by constructing a diffeomorphism
$\Lambda\to\mathsf{T}(H), \lambda\mapsto\varphi_\lambda$ from a
finite dimensional differentiable real manifold~$\Lambda$ of control
parameters~$\lambda$ and the optimized bounds for $e_0$ will be
\begin{equation}
  \sup_\Lambda\inf_\mathcal{Q}(E_{\varphi_\lambda})\leqslant e_0 \leqslant\inf_\Lambda\sup_\mathcal{Q}(E_{\varphi_\lambda})\;.
\end{equation}

\section{Inequalities in the differentiable case: old and new}\label{sec:oldnew}

\subsection{General considerations}

When $H$ is a local differential operator \textit{i.e.} involves a
finite number of derivatives in an appropriate representation (for
instance in position or in momentum representation), one can therefore
construct an algorithm that does not require any integration but
differential calculus only.  In an analytic and in a numerical
perspective, this may be a significant advantage on the perturbative
or variational methods even though it is immediate to see\footnote{
For each positive $\varphi$, it
  follows from 
\begin{multline*}
 \langle\varphi|H\varphi\rangle=\Re\big[\langle\varphi|H\varphi\rangle\big]=
  \int_\mathcal{Q}{\varphi}(q)\,\Re\big[H\varphi(q)\big]\,d\mu(q)\\
  =\int_\mathcal{Q}\varphi^2(q) E_{\varphi}(q)\,d\mu(q) \leqslant
  \langle\varphi|\varphi\rangle \sup_\mathcal{Q}
  \big(E_{\varphi}(q)\big).
\end{multline*}
} that the upper bound given by
\eqref{eq:inequalities} is always larger than
$\langle\varphi|H\varphi\rangle/\langle\varphi|\varphi\rangle$.  The
global analysis appears only through the determination of the
singularities and the \emph{absolute} extrema of the local energy that
may have bifurcated when the control parameter  $\lambda$ varies smoothly.  For a
Schr\"odinger operator, the possible singularities of the potential
on~$\bar{\mathcal{Q}}$ like a Coulombian divergence or an unbounded
behaviour at infinite distances may furnish a strong guideline for
constructing relevant test functions (see \S~\ref{subsec:HB} below).
I have given in \cite{Mouchet05a}, some heuristic and numerical
arguments to show how this strategy can be fruitful. In the present
article, the main focus will concern rigorous results and I will
explain how some of them can be obtained with great simplicity even
for systems as complex as those involved in the many-body problem.

\subsection{Case (i) Barta's inequalities} They immediately appear as a particular
case of theorem~\ref{th:main}:

\begin{theorem}[Barta 1937]\label{th:Barta} Let $\mathcal{Q}$ be a connected bounded 
  Riemannian manifold endowed with the metric~$g$ and $H$, the
  opposite of the Laplacian~$\Delta_g$ acting on the functions in
  $\mathsf{L}^2(\mathcal{Q},\mu)$ that satisfy Dirichlet boundary
  conditions on the boundary $\partial\mathcal{Q}$.  The lowest
  eigenvalue~$e_0$ of $H$ is such that, for all
  positive~$\varphi\in\mathcal{C}^2(\mathcal{Q})$,
\begin{equation}
  \inf_\mathcal{Q}\left(-\frac{\Delta_g\varphi}{\varphi}\right)\leqslant e_0 \leqslant\sup_\mathcal{Q}\left(-\frac{\Delta_g\varphi}{\varphi}\right)\;.
\end{equation}
\end{theorem}

\textit{Proof: } It follows directly from theorem~\ref{th:main}, with
the local energy given by $E_\varphi=-\Delta_g\varphi/\varphi$:
The spectrum of~$H$ is discrete and the Krein-Rutman theorem assures
that hypothesis~\ref{hyp:sympos} is fulfilled for $e_0$ being the
lowest (and simple) eigenvalue.  \qed

\begin{remark} The Dirichlet boundary
  conditions are not essential and can be replaced by any other type
  of boundary conditions provided that hypothesis~\ref{hyp:sympos}
  remains fulfilled. However, as explained in
  \S~\ref{subsec:strategy}, for obtaining interesting bounds on $e_0$
  extending $\mathsf{T}(H)$ to test functions that do not fulfill the
  boundary conditions (as proposed in \cite{Duffin47a}) seems not
  appropriate.
\end{remark}

\subsection{Case (ii) Duffin-Barnsley-Thirring inequalities}
Extensions of Barta's inequalities for Schr\"odinger operators have
been obtained partially by Duffin \cite{Duffin47a} and Barnsley
\cite{Barnsley78a} (for the lower bound only) and completely (lower
and upper bound) by Thirring \cite{Thirring79b} using Kato's
inequalities (see also \cite{Schmutz85a}).
\begin{theorem}[Duffin 1947, Barnsley 1978, Thirring 1979]\  \newline Let $H=-\Delta+V$ be
  a  Schr\"odinger operator acting on $\mathsf{L}^2(\mathbb{R}^d)$
  having an eigenvalue below the essential spectrum. Then the lowest
  eigenvalue~$e_0$ of $H$ is such that for any strictly positive
  $\varphi\in\mathsf{D}(H)$,
\begin{equation}
  \inf_{\mathbb{R}^d}\left(V-\frac{\Delta\varphi}{\varphi}\right)\leqslant e_0 
\leqslant\sup_{\mathbb{R}^d}\left(V-\frac{\Delta\varphi}{\varphi}\right)\;.
\end{equation}
\end{theorem}
The proof is similar to the one presented above with the local energy
being now $E_\varphi=V-\Delta\varphi/\varphi$. This argument has the advantage
on the existing ones that it does not involve the specific properties
of the Laplacian and can be immediately transposed to the larger class
of the differential operators (not necessarily of second-order) that
fulfill hypothesis~\ref{hyp:sympos}.

\subsection{Case (iii) Magnetic Schr\"odinger operators} \label{subsec:HB}
In the presence of a magnetic field, Schr\"odinger operators take the
form $H=(\imat\partial_q+A(q))^2+V(q)$ with $A:\mathcal{Q}\to\RR^d$
being a smooth magnetic potential vector and $V:\mathcal{Q}\to\RR$ a
smooth scalar potential.  The Krein-Rutman theorem may not apply
whereas there still exists a non-negative real eigenfunction ($e_0$
may be not simple nor the lowest eigenvalue)\cite{Helffer+99a}.

In the particular case of the Hydrogen atom in a constant and uniform
magnetic field, hypothesis~\ref{hyp:sympos} is fulfilled for all
values of the magnetic field \cite{Avron+77a,Avron+78a} and indeed
concerns the lowest eigenvalue. Therefore theorem~\ref{th:main}
applies and furnishes relevant analytical bounds that can be improved
numerically as shown in \cite{Mouchet05a}.
\begin{proposition} The smallest eigenvalue~$e_0$ of the (3d-)Zeeman Hamiltonian
\begin{equation}\label{eq:HamHB}
  H=\frac{1}{2}\left(-\imat\vec\nabla+\frac{1}{2}\vec{r}\times\vec{B}\right)^2-\frac{1}{r}
\end{equation}
is such that 
\begin{equation}
  \forall B\geqslant0,\quad e_0\leqslant-1/2+B/2.
\end{equation}
\end{proposition}

\textit{Proof: } In cylindrical coordinates ($\rho$, $\theta$, $z$)
where the magnetic field is~$\vec{B}=B\vec{u}_z$, the test function of
the form~$\varphi=\exp(-\sqrt{\rho^2+z^2}-B\rho^2/4)$ is constructed,
according to the strategy explained in \S~\ref{subsec:strategy}, in
order to respect the rotational invariance of the ground-state
\cite{Avron+77a} and eliminate both singularities at $r\to0$ and
at~$\rho\to\infty$. Indeed such a choice leads straightforwardly to
the bounded local energy:
\begin{equation}\label{eq:ElocHB}
   E_\varphi= -1/2+B/2-\frac{\rho^2B}{2\sqrt{\rho^2+z^2}}
\end{equation}
The upper bound follows.    \qed

\section{Application to the many-body problem}\label{sec:manybody}

\subsection{Expression of the local energy in terms of two-body functions}
We will consider in this section a $N$-body non-relativistic bosonic system in $d$
dimensions; $(d,N)\in(\NN\backslash\{0,1\})^2$; whose Hamiltonian is
given by:
\begin{equation}\label{eq:HNbody}
  \tilde{H}=\sum_{i=0}^{N-1}-\frac{1}{2m_i}\Delta_i+V(\vec{r}_0,\dots,\vec{r}_{N-1})
\end{equation} 
acting on $\mathsf{L}^2(\RR^{Nd})$, endowed with the canonical
Lebesgue measure, and where $\forall i\in\{0,\dots,N-1\}$,
$\vec{r}_i\in\RR^d$.  $\Delta_i$ is the Laplacian in the $\vec{r}_i$
variables and~$m_i\in\RR^+\backslash\{0\}$ the mass of the
$i^\text{th}$ particle. The spinless bosons interact only by the
two-body radial potentials~$v_{ij}=v_{ji}:\RR\to\RR$ \textit{i.e.} $V$
is given by
 \begin{equation}\label{eq:VNbody}
  V=\sum^{N-1}_{\substack{i,j=0\\i<j}}v_{ij}(r_{ij})
\end{equation}
where $r_{ij}=r_{ji}=||\vec{r}_j-\vec{r}_i||$.  Once the center of
mass is removed, the Hamiltonian~$\tilde{H}$ leads to a reduced
Hamiltonian~$H$ acting on $\mathsf{L}^2(\RR^{(N-1)d})$ (see for
instance \S~XI.5 of~\cite{Reed/Simon78b}) and we will suppose in the
following that $H$ has at least one eigenvector\footnote{Physically,
  this can be achieved with a confining external potential (a ``trap''
  is currently used in experiments involving cold atoms). Formally,
  this can be obtained in the limit of one mass, say $m_0$, being much
  larger than the others. The external potential appears to be the
  $v_{0,i}$'s, created by such an infinitely massive motionless
  device. It will trap the remaining~$N-1$ particles in some bounded
  states if the $v_{0,i}$'s increase sufficiently rapidly with the
  $r_{0,i}$'s.}. Therefore, hypothesis~\ref{hyp:sympos} is fulfilled
for $e_0$ being the lowest (and simple) eigenvalue.  A natural choice
for test functions is to consider factorized ones of the form:
\begin{equation}\label{eq:phiprod}
  \varphi(\vec{r}_0,\dots,\vec{r}_{N-1})=\prod^{N-1}_{\substack{i,j=0\\i<j}}\phi_{ij}(r_{ij})
\end{equation}
where $\phi_{ij}\in\mathsf{L}^2(\RR^+)$ and~$\phi_{ji}=\phi_{ij}>0$.
One can check easily that the total momentum of such a test function
vanishes.  The corresponding local energy~\eqref{def:localenergy} is
given by
\begin{equation}\label{eq:ElocNbody}\begin{split}
  E_\varphi(q_N)=\sum^{N-1}_{\substack{i,j=0\\i<j}}\left[
  \frac{-1}{2m_{ij}\phi_{ij}(r_{ij})}
   \left( \phi_{ij}''(r_{ij})+
  \frac{d-1}{r_{ij}}\,\phi_{ij}'(r_{ij})\right)+v_{ij}(r_{ij})\right]\\
  -\sum_{(\widehat{j,i,k})}\frac{1}{m_i}S'_{ij}(r_{ij})S'_{ik}(r_{ik})\cos(\widehat{j,i,k})
 \end{split}
\end{equation}
where $q_N\in\mathcal{Q}_N=\RR^{(N-1)d}$ stands for the $(N-1)d$
relative coordinates
$(\vec{r}_1-\vec{r}_0,\dots,\vec{r}_{N-1}-\vec{r}_0)$, $m_{ij}$ for
the reduced masses $m_im_j/(m_i+m_j)$, $S'_{ij}$ is the derivative of
$S_{ij}=\ln(\phi_{ij})$.  The last sum involves all the
$N(N-1)(N-2)/2$ angles $(\widehat{j,i,k})$ between
$\vec{r}_j-\vec{r}_i$ and $\vec{r}_k-\vec{r}_i$
 that can be formed with all
the triangles made of three distinct particles.

Whenever each $v_{ij}$ allows for a two-body bounded state, we can choose
$\phi_{ij}$ to be the eigenvector of $-(2m_{ij})^{-1}\Delta +v_{ij}$
having the smallest eigenvalue $\epsilon^{(2)}_{ij}$. Moreover, if
$v_{ij}(r)$ is bounded when $r\to\infty$, from an elementary semiclassical analysis (see
for instance \cite{Maslov/Fedoriuk81a}) it follows that~$S'_{ij}$ is
also bounded since asymptotically we
have~$S'_{ij}(r)\sim_{r\to\infty}-\sqrt{2m_{ij}[v_{ij}(r)-\epsilon^{(2)}_{ij}]}$.
It follows that the local energy is also bounded and finite lower and
upper bounds on~$e_0$ can be found. For instance, directly from
expression~\eqref{eq:ElocNbody}, we have

\begin{proposition} If, for all $(i,j)\in\{0,\dots,N-1\}^2$, $i\neq j$, $v_{ij}(r)$ is bounded
  when $r\to\infty$ and $-(2m_{ij})^{-1}\Delta +v_{ij}$ has a smallest
  eigenvalue~$\epsilon^{(2)}_{ij}$ obtained
  for~$\phi_{ij}=\exp{S_{ij}}$, then the smallest eigenvalue~$e_0$ of
  the N-body Hamiltonian in the center-of-mass frame is bounded by
\begin{equation}\label{eq:ineqNbody}
  \sum^{N-1}_{\substack{i,j=0\\i<j}}\epsilon^{(2)}_{ij}-\frac{s^2}{2m}\,N (N-1)(N-2) \leqslant
e_0 \leqslant \sum^{N-1}_{\substack{i,j=0\\i<j}}\epsilon^{(2)}_{ij}+\frac{s^2}{2m}\,N (N-1)(N-2). 
\end{equation}
where  $m=\min_{i}m_i$ and $s=\max_{i,j}\sup_{\RR^+}|S'_{ij}|$.
\end{proposition}

For potentials that are relevant in physics (see for instance the
effective power-law potentials of the form
$v_{ij}(r)=\mathrm{sign}(\beta)r^\beta$; $\beta\in\RR$; between
massive quarks as studied in \cite{Benslama+98a}), the analytic form
of the two-body eigenvector is not known in general and some numerical
computations are required to obtain the absolute maximum and minimum
of the local energy~\eqref{eq:ElocNbody}.

\subsection{The local energy for a general Coulombian problem} When
$N$ is large, the estimation~\eqref{eq:ineqNbody} is quite rough, in
particular
it does not take into account the constraints between the several angles.
More precise results are obtained for the Coulombian problem where
$v_{ij}(r)=e_{ij}/r$ with $e_{ij}\in\RR$ and~$d>1$. In that case,
provided a bounded state exists and that we keep the test
function~\eqref{eq:phiprod} in~$\mathsf{L}^2(\RR^{(N-1)d})$, we choose
a constant derivative~$S'_{ij}=2e_{ij}m_{ij}/(d-1)$ in order to get
rid of the Coulombian singularities of~$V$.  We obtain a bounded local
energy given by
\begin{equation}\label{eq:ElocNbodyCoulomb}
  E_\varphi(q_N)=\sum^{N-1}_{\substack{i,j=0\\i<j}}-\frac{2m_{ij}e_{ij}^2}{(d-1)^2}-
\frac{4}{(d-1)^2}\sum_{(\widehat{j,i,k})}\frac{m_{ij}m_{ik}e_{ij}e_{ik}}{m_i}\cos(\widehat{j,i,k})\;.
\end{equation}

\subsection{Identical purely attractive Coulombian particles} The case
where all the $N$ particles 
are identical and attract each other, \textit{i.e.} when
$\forall(i,j)\in\{0,\dots,N-1\}^2$, $i\neq j$, $m_i=1$ and $e_{ij}=-1$,
has been extensively studied in the literature, in particular the
asymptotic behavior of~$e_0$ with large~$N$ may have some dramatic
consequences on the thermodynamical limit
\cite{Fisher/Ruelle66a,Lenard/Dyson67a,LevyLeblond69a}.  The local
energy method allow to obtain in a much simpler way energy bounds that
are comparable to those already obtained by other methods.  Actually,
\eqref{eq:ElocNbodyCoulomb} simplifies to
\begin{equation}\label{eq:ElocNbodyidCoulomb}
  E_\varphi(q_N)=-\frac{1}{(d-1)^2}\left[\frac{1}{2}\,N(N-1)+F_N(q_N)\right]
\end{equation} 
with $F_2\equiv0$ and for~$N\geqslant3$
\begin{equation}\label{def:FN}
  F_N(q_N)=\sum_{(\widehat{j,i,k})}\cos(\widehat{j,i,k}).
\end{equation} 
The angular function~$F_N$ depends only on the geometrical
configuration of the $N$ vertices $(\vec{r}_0,\dots,\vec{r}_{N-1})$
\textit{i.e.}  it is invariant under the group of Euclidean isometries
and the scale invariance of the Coulombian interaction makes it
invariant under dilations as well.
\begin{lemma}\label{lem:F3}
  $\inf_{\mathcal{Q}_3}F_3=1$ is obtained when the three points are
  aligned.  $\sup_{\mathcal{Q}_3}F_3=3/2$ is obtained when the three
  points make an equilateral triangle.
\end{lemma}
\textit{Proof: }The extrema of $F_3$ correspond to the extrema of the
function defined by
$(\theta_0,\theta_1,\theta_2)\mapsto\cos\theta_0+\cos\theta_1+\cos\theta_2$
under the constraint~$\theta_0+\theta_1+\theta_2=\pi$ for
$(\theta_0,\theta_1,\theta_2)$ being the three
angles~$(\widehat{1,0,2})$, $(\widehat{0,1,2})$, $(\widehat{0,2,1})$
respectively.  The Lagrange multiplier method leads to the
determination of the extrema of the function $[0;\pi]^3\to\RR$ defined
by
$(\theta_0,\theta_1,\theta_2)\mapsto\cos\theta_0+\cos\theta_1+\cos\theta_2+l(\theta_0+\theta_1+\theta_2-\pi)$;
$l\in\RR$. We immediately obtain that the extremal points are located
at~$(\theta_0,\theta_1,\theta_2)=(\pi/3,\pi/3,\pi/3)$ and
$(\theta_0,\theta_1,\theta_2)=(0,0,\pi)$ together with the solutions
that are obtained by circular permutations. It is easy to check that
the first solution provides an absolute maximum for~$F_3$ and the
second ones an absolute minimum. \qed

An immediate consequence of the preceding lemma is
\begin{proposition}The lowest energy $e_0$ of $N\geqslant2$ 
identical attractive Coulombian spinless particles in $d>1$ dimensions
  is such that
\begin{equation}\label{eq:upperENbodycoulombian}
  e_0\leqslant-\frac{1}{(d-1)^2}\frac{1}{6}\,N(N-1)(N+1)
\end{equation} 
when the individual masses equal to unity and the attractive potential is~$-1/r$. 
\end{proposition}
\textit{Proof: }The sum on the angles that defines $F_N$ in
\eqref{def:FN} can be written as a sum of $N(N-1)(N-2)/6$ $F_3$-terms
 calculated for all the triangles that belong to the $N$-uplet made of the $N$ vertices.
Since, from Lemma~\ref{lem:F3} the absolute minimum of $F_3$ is
obtained for a flat configuration, when all the~$N$ points are aligned
all the~$F_3$-terms reach their absolute minimum simultaneously and
the absolute minimum of~$F_N$ is obtained. We have
\begin{lemma}
  $\inf_\mathcal{Q_N}F_N=\displaystyle\frac{1}{6}\,N(N-1)(N-2)$ is obtained when all the $N$ points are aligned.
\end{lemma}
The upper bound~\eqref{eq:upperENbodycoulombian} follows from
\eqref{eq:ElocNbodyidCoulomb}. \qed

More generally, for a given $N$-uplet (\textit{i.e.} a set of exactly
$N$ points), clustering the sum \eqref{def:FN} in $M$-uplets
($N\geqslant M\geqslant 3$) allows to find bounds on~$F_N$ from bounds
on~$F_M$. Indeed, we can write
\begin{equation}\label{eq:clustering}
  F_N(q_N)=\sum_{q_M}\frac{(M-3)!\,(N-M)!}{(N-3)!}F_M(q_M)
\end{equation}
where the sum is taken on all the $M$-uplets, labeled by the
coordinates~$q_M$, that belong to the given~$N$-uplet. This sum
involves exactly $N!/M!/(N-M)!$ terms and  we have, therefore,
 \begin{lemma}$\forall(N,M)\in\NN^2\ \text{such that}\  N\geqslant M\geqslant3$,
\begin{equation}
  \sup_{\mathcal{Q}_N}F_N\leqslant\frac{N(N-1)(N-2)}{M(M-1)(M-2)}\,\sup_{\mathcal{Q}_M}F_M\;.
\end{equation}
\end{lemma}
For a given~$N$, $\sup_{\mathcal{Q}_N}F_N$ is not known exactly but the ordered sequence 
\begin{multline}
\frac{\sup_{\mathcal{Q}_N}F_N}{N(N-1)(N-2)}\leqslant\frac{\sup_{\mathcal{Q}_{N-1}}F_{N-1}}{(N-1)(N-2)(N-3)}
\leqslant\cdots\\
\cdots\leqslant\frac{\sup_{\mathcal{Q}_M}F_M}{M(M-1)(M-2)}\leqslant\cdots\leqslant
\frac{\sup_{\mathcal{Q}_3}F_3}{3.2.1}=\frac{1}{4} 
\end{multline}
shows that in order to improve the lower bounds
on~\eqref{eq:ElocNbodyidCoulomb}, we must try to find
$\sup_{\mathcal{Q}_M}F_M$ with $M$ being the largest as possible.
However, when considering identity~\eqref{eq:clustering} for~$N=4$
and~$M=3$ together with Lemma~\ref{lem:F3} we have
\begin{lemma}\label{lem:F4}
  $\sup_{\mathcal{Q}_4}F_4=6$ is obtained when the four points make a regular tetrahedron.
\end{lemma}
Then
$\sup_{\mathcal{Q}_4}F_4/(4.3.2)=\sup_{\mathcal{Q}_3}F_3/(3.2.1)=1/4$
and no better estimate is obtained when considering~$M=4$ rather
than~$M=3$.  Numerical investigations lead to the following
conjectures:
\begin{conjecture}[$C_5$] When~$d=3$,
$\displaystyle\sup_{\mathcal{Q}_5}F_5
=\frac{9}{2}+\frac{6(h_0+1)}{\sqrt{h_0^2+\frac{1}{3}}}-\frac{1}{h_0^2+\frac{1}{3}}\simeq14.591594$
with 
\begin{multline}\label{eq:h0}
  6h_0=1+\sqrt {-1+\sqrt [3]{7+4\,\sqrt {3}}+{\frac {1}{\sqrt [3]{7+4
\,\sqrt {3}}}}}\ +\ \\
\sqrt {-2-\sqrt [3]{7+4\,\sqrt {3}}-{\frac {1}{
\sqrt [3]{7+4\,\sqrt {3}}}}+8\,{\frac {1}{\sqrt {-1+(7+4\,
\sqrt {3})^{1/3}+(7+4\,\sqrt {3})^{-1/3} }}}}
\end{multline}
is obtained when the five points make two mirror-symmetric
tetrahedrons sharing one common equilateral basis, their other faces
being 6 isosceles identical triangles.
\end{conjecture}
\begin{remark} The only free parameter of the specific configuration 
  can be chosen to be the height~$h$ of one tetrahedron (the length of
  the edges of the common equilateral basis being fixed to one).  The
  maximum of $h\mapsto F_5$ is reached for~$h_0$ being the greatest
  solution of~$9\,h^4-6\,{h}^{3}+3\,{h}^{2}-2\,h+1/3$ that is
  precisely given by~\eqref{eq:h0}.
\end{remark}
\begin{remark} The  pyramidal configurations with a squared basis 
  leads to a local maximum that gives $F_5=15/2+5\sqrt{2}\simeq14.57$.
\end{remark}
\begin{conjecture}[$C_6$]When~$d=3$,
  $\sup_{\mathcal{Q}_6}F_6=12(1+\sqrt{2})$ is obtained when the 6
  points make a regular octahedron.
\end{conjecture}
\begin{conjecture}[$C_8$]When~$d=3$, $\sup_{\mathcal{Q}_8}F_8={\scriptstyle16}[
  \frac{4}{5}+\frac{1}{\sqrt{2}}+\frac{1}{\sqrt{5}}+
  \frac{4(1+\sqrt{2})}{\sqrt{5}\sqrt{5+4\sqrt{2}}}+
  \frac{3+2\sqrt{2}}{\sqrt{5+4\sqrt{2}}}- \frac{1}{5+4\sqrt{2}}
  ]\simeq79.501 $ is obtained when the 8 points make two identical
  squares (whose edges have length one) lying in two parallel planes
  separated by a distance $h=\sqrt{1+2\sqrt{2}}/2$. The axis joining
  the centers of the two squares is perpendicular to the squares and
  the two squares are twisted one from the other by a relative angle
  of $\pi/4$.
\end{conjecture}
\begin{remark} The cube corresponds to $F_8=8(3\sqrt{2}+\sqrt{3}+3/2+\sqrt{6})\simeq79.393$.
\end{remark}
\begin{conjecture}[$C_\infty$]When~$d=3$ and $N\to\infty$, the configuration that maximizes $F_N$ corresponds to
  $N$ points uniformly distributed on a sphere and
  $\sup_{\mathcal{Q}_N}F_N\sim\frac{2}{9}N^3+o(N^3)$.
\end{conjecture}
\begin{remark} The ambiguity of distributing $N$ points uniformly on a
sphere \cite[and references therein]{Saff/Kuijlaars97a} vanishes for large $N$
as far as a uniform density is obtained. Assuming such a uniform density, the continuous limit of~$F_N/N^3$
is a triple integral on the sphere than can be computed exactly to $2/9$.  
\end{remark}

The upper bound \eqref{eq:upperENbodycoulombian} for $d=3$ is slightly
above the one obtained 
in \cite[eq. (17 p. 807)]{LevyLeblond69a},
namely $-(5/8)^2 N(N-1)^2/8$, and the numerical estimate
$-.0542N(N-1)^2$ 
in \cite[eq. (16) p. 63]{Basdevant+90a}. For the lower bounds, from 
Lemmas~\ref{lem:F3} or~\ref{lem:F4}, we have obtained
 \begin{proposition}The lowest energy $e_0$ of $N\geqslant2$ 
identical attractive Coulombian spinless
   particles in $d>1$ dimensions is such that
\begin{equation}\label{eq:lowerENbodycoulombian}
  -\frac{1}{(d-1)^2}\frac{1}{4}\,N^2(N-1)\leqslant e_0
\end{equation} 
when the individual masses equal to unity and the attractive potential is~$-1/r$. 
\end{proposition}
The same result has been obtained for~$d=3$ in \cite[eq. (11) p. 62]{Basdevant+90a} and is slightly 
better than
 $-N(N-1)^2/8$ given in \cite[eq. (13 p. 807)]{LevyLeblond69a}. For~$N\geqslant M\geqslant4$, this lower bound can be
improved to
\begin{equation}
   -\frac{1}{(d-1)^2}\,N(N-1)(\frac{1}{2}+\alpha_M(N-2))\leqslant e_0
\end{equation}
with
\begin{equation}  
\alpha_M=\frac{\sup_{\mathcal{Q}_M}F_M}{M(M-1)(M-2)}\leqslant\frac{1}{4}. 
\end{equation}
If conjecture $(C_5)$ [resp. $(C_6)$, $(C_8)$ and $(C_\infty)$  ] is
correct we get, when $d=3$, 
$\alpha_5\simeq.2432$ [resp. $\alpha_6\simeq.2414$, $\alpha_8\simeq.2366$ and $\alpha_\infty\sim2/9$ ]
and  the lower bounds are, therefore, improved.

Some numerical investigations, in particular a systematic comparison
with the lower bounds obtained with variational methods  in \cite{Benslama+98a} for $N=3$ and
$N=4$ Coulombian particles  will be given elsewhere \cite{Mouchet05c}.  

\section{Application to discrete Hamiltonians}\label{sec:discrete}
 
Where~$q\in\ZZ^d$; $d\in\NN$; the discretized analogue of a local
differential operator corresponds to a Hamiltonian that couples at
most a finite number of basis vectors (\textit{e.g.} the nearest neighbors on the
lattice~$\ZZ^d$). For instance, when~$d=1$, it can be seen as a
Hermitian band matrix (finite or infinite) of finite half-width in an
appropriate basis. Possibly with renumbering the~$q$'s, on
$\ell^2(\ZZ^d)$ $H$ has the form:
\begin{equation}\label{eq:discreteham}
  (H\varphi)_q=\sum_{\substack{\nu\in\ZZ^d\\|\nu|_\infty\leqslant N_b}}H_{q,q+\nu}\varphi_{q+\nu}\;,
\end{equation}
where $N_b\in\NN$, $\forall(q,q')\in\ZZ^d\times\ZZ^d$,
$H_{q',q}=\overline{H}_{q,q'}\in\CC$ and $|\nu|_\infty$ stands for
$\max(|\nu_1|,\dots,|\nu_d|)$. $\varphi$ will be taken as a discrete
set of real strictly positive numbers, and the local
energy~$E_\varphi(q)$ is computed, for a given $q$, with elementary
algebraic operations whose number is finite and all the smaller
than~$N_b$ is small: its value at a given $q$ depends on~$(2N_b)^d+1$
components of~$\varphi$ at most.  Under
hypothesis~\ref{hyp:sympos}\footnote{We have seen at the end of
  section~\ref{subsec:genineq} that the symmetry hypothesis can be
  relaxed.}, if, say, for a given test vector $\varphi$, the absolute
maximum of $E_\varphi$ occurs only at a unique finite~$q_m$, one can
immediately improve the upper bound by a finite amount, for instance
just by varying $\varphi_{q_m}$ only, until $E_\varphi(q_m)$ is not an
absolute maximum anymore. Only $(2N_b)^d+1$ values of the local energy
will be affected by the variation of just one component of $\varphi$.
One can see easily that this approach leads to a wide variety of
algorithms where a sequence of optimization steps is constructed; each
step involves a number of optimization parameters and functions that
is usually much smaller (of order $(2N_b)^d$ or less) than the
dimension of the original matrix.

Discrete Schr\"odinger operators are important particular cases of
Hamiltonians~\eqref{eq:discreteham} with $N_b=1$. They are relevant
models for the description of quantum (quasi-)particles evolving in
periodic crystals.  For $d=1$, they can be written as
\begin{equation}\label{eq:discreteschro}
  (H\varphi)_q=-\varphi_{q+1}-\varphi_{q-1}+V(q)\,\varphi_q
\end{equation}
where the potential~$V$ is a real bounded function on $\ZZ$. By
possibly subtracting a constant positive real number to $V$, the
bounded operator $-H$ can be made
positive and ergodic: First, for any positive and non-identically
vanishing $\varphi$ and $\varphi'$ in $\ell^2(\ZZ)$, $-H\varphi$
remains positive. Second,
$\langle\varphi'|(-H)^{|q'-q|}\varphi\rangle=\varphi'_{\smash[t]{q'}}\varphi\vphantom{\varphi'}_{q\vphantom{{q'}}}+(\text{positive
  terms})\neq0$ for any given pair of strictly positive components
$\varphi'_{\smash[t]{q'}}$ and $\varphi_{q}$ with $q'\neq q$ (In the marginal case where
$\varphi$ and $\varphi'$ both vanish everywhere but on the same point,
we have $\langle\varphi'|(-H)\varphi\rangle)\neq0$). 
Therefore theorem XIII.43 of
\cite{Reed/Simon78a} applies: if~$H$ has indeed one eigenvalue,
hypothesis~\ref{hyp:sympos} is fulfilled for $e_0$ being the smallest
eigenvalue of $H$.  In that case, inequalities~\eqref{eq:inequalities}
takes the form:

\begin{proposition} When the discrete Schr\"odinger
 operator~\eqref{eq:discreteschro} admits at least
  one eigenvalue and when $V$ is bounded, then the smallest
  eigenvalue~$e_0$ is such that, $\forall\varphi\in\ell^2(\ZZ)\ 
  \mathrm{such\ that}\ \varphi>0$,
\begin{equation}
\inf_{q\in\ZZ}\left(-\frac{\varphi_{q+1}+\varphi_{q-1}}{\varphi_q}+V(q)\right)
\leqslant e_0 \leqslant
\sup_{q\in\ZZ}\left(-\frac{\varphi_{q+1}+\varphi_{q-1}}{\varphi_q}+V(q)\right).
\end{equation}
\end{proposition}

When $V$ is actually a $N$-periodic real function,
the spectral problem (see \cite{Reed/Simon78a} for instance) leads 
to the search of complex series~$(u_q)_{q\in\ZZ}$ such that
\begin{equation}\label{eq:bandspectrum}
\forall\eta=(\eta_1,\eta_2)\in[0;1[^2\quad\left\{
\begin{aligned}&-u_{q+1}-u_{q-1}+V(q+\eta_2)\,u_q=e(\eta)u_q\;;\\
             &u_{q+N}=\EXP{\imat2\pi\eta_1}u_{q}\;.
\end{aligned}\right.
\end{equation}

The spectrum of~$H$ is the bounded
set~$\sigma(H)=\{e(\eta)|\eta\in[0;1[^2\}\subset\RR$.  It is
given by the reunion for all $\eta$'s of the $N$ eigenvalues of
finite $N\times N$ Hermitian matrices~$H^{(\eta)}$ obtained after
transforming \eqref{eq:bandspectrum} with the one-to-one
mapping~$u_{q}\mapsto u_{q}\exp(-\imat 2\pi q\eta_1/N)$. As far as
positive solutions of~\eqref{eq:bandspectrum} are concerned, we will
take~$\eta_1=0$ and will look for the smallest
eigenvalue~$e_0(\eta_2)$ of
\begin{equation}\label{eq:Htheta2}
  H^{(0,\eta_2)}=\begin{pmatrix}
V(\eta_2)&-1&0&\cdots&0&-1\\
-1&V(1+\eta_2)&-1&0&\cdots&0\\
0&-1&V(2+\eta_2)&-1&\hdotsfor{2}\\
0&\hdotsfor{5}\\
-1&0&\cdots&0&-1&V(N-1+\eta_2)\\
\end{pmatrix}
\end{equation}
\begin{remark}
  The (rational) Harper model \cite{Harper55a,Harper55b} (also called
  the almost Mathieu equation~\cite{Bellissard/Simon82a}) corresponds
  to $V(q)=-V_0\cos(2\pi q M/N)$ where $V_0>0$, $(M,N)$ being strictly
  positive coprimes integers. For a given $N$ and~$M$, $\sigma(H)$
  appears to be made of~$N$ bands.  The union of these bands for each
  rational number~$M/N$ between 0 and 1 produces the so-called
  Hof\-stad\-ter butterfly~\cite{Hofstadter76a}.
\end{remark}
We are therefore able to produce two non-trivial bounds on the lowest
eigenvalue~$e_0(\eta_2)$ without any diagonalization:

\begin{proposition} When~$V$ is $N$-periodic,
  $\forall\eta_2\in[0;1[$, 
the smallest eigenvalue~$e_0(\eta_2)$ 
of~\eqref{eq:Htheta2} is such that 
 $\forall\varphi\in(\RR^+\backslash\{0\})^N$,
\begin{subequations}\label{subeq:inequalitiesdiscrete}
\begin{equation}
\min_{q\in\{0,\dots,N-1\}}\left(-\frac{\varphi_{q+1}+\varphi_{q-1}}{\varphi_q}+V(q+\eta_2)\right)
\leqslant e_0(\eta_2)
\end{equation}
and
\begin{equation}
e_0(\eta_2) \leqslant
\max_{q\in\{0,\dots,N-1\}}\left(-\frac{\varphi_{q+1}+\varphi_{q-1}}{\varphi_q}+V(q+\eta_2)\right)
\end{equation}
\end{subequations}
(the indices labeling the components of $\varphi$ are taken modulo $N$).
\end{proposition}

Therefore, we can bound the
bottom of the Hofstadter butterfly with the help of any test function:
\begin{corollary}For the rational Harper Hamiltonian
\begin{equation}
  H\varphi_q=-\varphi_{q+1}-\varphi_{q-1}-V_0\cos\left(2\pi q \frac{M}{N}\right)\,\varphi_q\;,
\end{equation}
we have $\forall\varphi\in(\RR^+\backslash\{0\})^N$
\begin{subequations}\label{subeq:inequalitiesharper}
\begin{equation}
\min_{q\in\{0,\dots,N-1\}}\left[-\frac{\varphi_{q+1}+\varphi_{q-1}}{\varphi_q}-V_0\cos\left(2\pi q \frac{M}{N}\right)\right]
\leqslant \inf\left[\sigma(H)\right]
\end{equation}
and
\begin{equation}
\inf\left[\sigma(H)\right] \leqslant
\max_{q\in\{0,\dots,N-1\}}\left[-\frac{\varphi_{q+1}+\varphi_{q-1}}{\varphi_q}-V_0\cos\left(2\pi q\frac{M}{N} \right)\right]\;.
\end{equation}
\end{subequations}
\end{corollary} 

\textit{Proof: }In the Harper model, for each rational number~$M/N$,
the lowest eigenvalue is obtained for $\eta=(0,0)$. It is a direct
application of \cite[Theorem~XIII.89~(e)]{Reed/Simon78a} and thus
inequalities~\eqref{subeq:inequalitiesharper} follow directly
from~\eqref{subeq:inequalitiesdiscrete}.\qed

Choosing for $\varphi$, at first guess, semiclassical approximations
(\textit{i.e.} corresponding to large~$N$) constructed from Mathieu
functions are therefore expected to provide numerical reasonable
bounds.

\section{Conclusion}

It has been shown on various examples how theorem~\ref{th:main} can be
used to obtain rigorous estimates on the principal value of any
symmetric operator. Its simplicity, its low cost in computations and
its wide domain of applications make the method presented in this
article a powerful tool for controlling bounds. In many situations, it
provides non-trivial complementary information to those obtained by
traditional or more system-dependent methods. Unfortunately, I was not
able to extend the method to fermionic systems (see for
instance~\cite{Sigal95a} and references therein) where the spatial
wavefunction of the ground-state has generically non-trivial
nodes~\cite{Ceperley91a} that cannot be known \textit{a priori} even
with some considerations on symmetries.

I have also presented some clues for further developments of
optimization algorithms. However it remains an open question whether
such algorithms really bear the potential of an efficient treatment
and will overcome the possible difficulties one may face in realistic
problems.

\acknowledgement{} I thank Laurent V\'eron for a critical reading of
the first proofs of this manuscript, Evans Harrell for fruitful
discussions, Sa\"{\i}d Ilias for letting me know useful references on
the generalizations of Barta's inequalities, Emmanuel Lesigne for
attracting my attention to reference~\cite{Saff/Kuijlaars97a} and
Michel Caffarel for reference~\cite{Ceperley91a}. I also gratefully
acknowledge the Laboratoire Kastler Brossel for its warm hospitality.


\newcommand{\etalchar}[1]{$^{#1}$}

\end{document}